\documentclass[letterpaper, 10pt]{article}      

\usepackage{graphicx}    
\usepackage{amsmath}
\usepackage{amssymb}
\usepackage{gensymb}
\usepackage{url}
\usepackage[english]{babel}
\usepackage{cite}
\usepackage{color,soul}
\usepackage{setspace}

\newtheorem{theorem}{Theorem}
\newcommand{\done}{\hspace*{\fill} $\Box$}

\title{\LARGE \bf Flight Management System for Hydrogen-Powered Aircraft in Cruise}

\author{Maxim Kaptsov and Luis Rodrigues\\
	\thanks{$^{1}$PhD student
		{\tt\small maxim.kaptsov@concordia.ca}}%
	\thanks{$^{2}$Professor
		{\tt\small luisrod@encs.concordia.ca}}%
	{\normalsize\itshape
		Department of Electrical and Computer Engineering}\\
	{\normalsize\itshape
		Concordia University, Montreal, Quebec, H3G 2W1, Canada}\\
}

\onecolumn
\begin{document}
\maketitle
\thispagestyle{empty}
\pagestyle{empty}

\begin{abstract}
The minimization of the Direct Operating Cost (DOC) for hydrogen-powered aircraft is formulated in this paper as an optimal control problem and is solved based on Pontryagin's minimum principle.
As a consequence, the optimum cruise flight speed is determined assuming cruising at a constant altitude. The optimization criterion corresponds to the minimization of a functional representing the trade-off between the cost of hydrogen fuel and time-dependent costs, which are related by a parameter denoted by cost index. The value of this parameter is introduced by a pilot into the Flight Management System (FMS) of the aircraft. The HY4 aircraft model is used to obtain numerical results for the proposed methodology.
\end{abstract}

\noindent{\bf Keywords:} Flight management systems, hydrogen fuel cell aircraft, optimal control, Pontryagin's principle

\section{Introduction}
Each new generation of commercial aircraft is designed to be more environmentally friendly than the previous one. However, total aircraft emissions keep increasing because the air traffic grows at a high rate. For example, the International Air Transport Association (IATA) predicted that the number of air passengers will double by 2036 \cite{IATA}. As a consequence, numerous startup projects have emerged recently (some of which are backed by manufacturers such as the Boeing Company and Airbus) aimed at developing emission-free passenger aircraft. On a quest for greener skies, several concepts have been proposed and studied in the last few years, including all-electric battery-powered aircraft (such as NASA X-57 Maxwell) and hybrid electric aircraft (Airbus E-Fan X). Another prominent example are airplanes powered by hydrogen fuel cells. 

A fuel cell is essentially an electrochemical device which converts hydrogen into electrical energy. The only exhaust products it generates are heat and water vapor. Other than being nearly emission-free, a fuel cell aircraft produces zero noise which is especially important during takeoff and landing in urban centers.
Due to the various technical challenges (some of which are related to water and heat management), it is unlikely that this technology would find its way into large commerical aircraft in the next decade. However, an airplane powered by hydrogen fuel cells is certainly a viable solution for aerobatic aircraft and the emerging business of urban mobility with so-called air taxis.
In fact, future air taxis have significant potential but some research issues must be solved first. One important research aspect is to be able to reduce the noise of such aircraft because the operation of air taxis will be intended for city centers. Hydrogen-powered aircraft then become a potential interesting solution if trajectories that minimize the consumed energy can be developed.
This optmization is the main theme of the paper. 

Some research on application of hydrogen cell technology in aircraft appeared during the decade 2000--2010. The authors of \cite{Bataller} studied an architecture of the power source which included a polymer electrolyte membrane fuel cell to power the Boeing FCD airplane. In \cite{Furrutter}, the authors proposed the design of a power plant based on hydrogen fuel cells to propel a small unmanned aerial vehicle. The design of a 2-seater electrical airplane powered by fuel cells was presented in \cite{Romeo}. In \cite{Jenal}, the authors studied performance characteristics and effectiveness of the proton exchange membrane fuel cell propulsion system for small aircraft. The feasibility of the fuel cell stack as an auxillary power unit in more electric aircraft was discussed in \cite{Dai}, while \cite{Kaptsov1, Kaptsov2} obtained the optimal cruising speed of a battery-powered airplane. 
More recently, research on hybrid aircraft that use hydrogen fuel cells has been reported in reference \cite{Salehpour}. The authors concluded that their optimal solution for fuel cell power can improve the aircraft performance measured in terms of fuel consumption and emission of polluted gas to the atmosphere. A recent study in reference \cite{Boll} focuses on the use of liquid hydrogen for cryogenic electric propulsion of hybrid electric aircraft.
A comprehensive study and comparison of electric, hybrid, and turboelectric fixed-wing aircraft can be found in \cite{Brelje}.

To the best of our knowledge, the trade-off optimality of energy use versus flight time has not been addressed in the open literature of hydrogen-powered vehicles.
Determining the optimal speed in the presence of a trade-off between fuel costs and travel time-related costs defines the problem that is formulated and solved in this work.
The rest of the paper is organized as follows: Section II presents the dynamic model of a fuel cell aircraft cruising at a constant altitude. In Section III, this model is used to formulate an optimal control problem. Pontryagin's minimum principle is then employed to solve the optimal control problem. The proposed method is applied in Section IV to the HY4 airplane.

\section{Dynamic model of a fuel cell aircraft}
\subsection{Longitudinal flight dynamics}
The nonlinear longitudinal dynamic model of an aircraft is given by the following differential equations [\citen{Anderson}]
\begin{equation}
\label{Dynamics}
\begin{split}
\dot{x}&=v\cos \gamma\\
\dot{h}&=v\sin \gamma\\
\dot{v}&=\left( \frac{g}{W} \right)\left( T\cos \alpha-D-W\sin \gamma \right)\\
\dot{\gamma}&=\left( \frac{g}{Wv} \right)\left( T\sin \alpha +L-W\cos \gamma \right)\\
\end{split}
\end{equation}
where $x$ is the horizontal position, $h$ is the altitude, $v$ is the true airspeed, $\gamma$ is the flight path angle, $\alpha$ is the angle of attack, $T$ is the thrust, $W$ is the weight, $L$ is the lift,  and $g$ is the gravitational acceleration.
Without loss of generality, the initial horizontal position is assumed to be $x(0)=0$ and the final position is assumed to be gven as $x(t_f)=x_d$. 

An equation for drag was developed in [\citen{Anderson}] for flight below drag divergence Mach number and is given by
\begin{equation}\label{Drag}
D=\frac{1}{2}C_{D,0}\rho Sv^2+\frac{2KW^2}{\rho Sv^2}
\end{equation}
where $C_{D,0}$ is the zero-lift drag coefficient and $K$ is the induced drag coefficient.
The following assumptions will be used to simplify the model (\ref{Dynamics})-(\ref{Drag})

{\bf Aircraft Assumptions:}
\begin{enumerate}
	\item The aircraft cruises at a constant predetermined altitude. Therefore, in cruise $\gamma=\dot{\gamma}=\dot{h}=0$.
	\item The angle of attack $\alpha$ is small. Therefore, $\cos\alpha\approx1$, $\sin\alpha\approx\alpha$.
	\item The component of the thrust perpendicular to the velocity vector is small compared to lift L and weight W.
	\item The flight Mach number is less than the drag divergence Mach number. 
	\item The aircraft is in steady cruise and the objective is to determine the flight speed $v$.
\end{enumerate}
Using these assumptions, the model (\ref{Dynamics}) can be simplified to
\begin{equation}
\label{Dynamics_reduced}
\begin{split}
\dot{x}&=v\\
\dot h&=0\\
L&=W\\
T&=D\\
\end{split}
\end{equation}

In this dynamic model $v$ acts as a control input.
\subsection{Fuel cell dynamics}
\begin{figure}[t]
	\centering
	\includegraphics[width=1\linewidth]{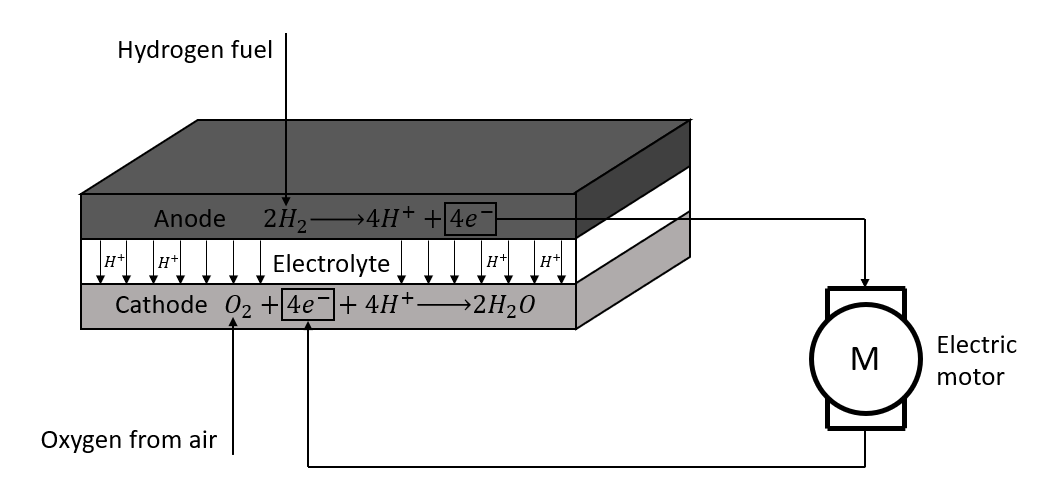}
	\caption{A block diagram of a fuel cell (adapted from \cite{Larminie})}
	\label{fig:screenshot005}
\end{figure}

A fuel cell converts chemical energy into electrical energy through what is called a reduction-oxidation (redox) chemical reaction.
For a hydrogen-oxygen fuel cell there is an anode and a cathode, with oxidation of the hydrogen occurring at the anode and producing both electrons and protons.
A block diagram of a fuel cell is shown in Figure \ref{fig:screenshot005} where the reactions of oxidation at the anode and reduction at the cathode are detailed. 
Electrons released in the chemical reaction go through an electrical circuit from the anode to the cathode.
The protons move through an electrolyte to the cathode in the opposite direction of the electrons.
At the cathode the electrons and protons combine with each other and the oxygen to produce water.
The operation of a hydrogen fuel cell can be represented by the following simple chemical equation combining the cathode and anode reduction-oxidation reactions, which corresponds to the addition of the two equations in Figure \ref{fig:screenshot005}
\begin{equation}\label{chemical_reaction}
2H_2 + O_2 \longrightarrow 2H_2O
\end{equation}

Note that two hydrogen molecules are required for each oxygen molecule. 
By Faraday's laws, the charge produced in a single fuel cell is directly proportional to the amount of hydrogen consumed and given by
\begin{equation}\label{single_cell_charge}
Q = 2F\frac{m_{H}}{M_{H}}
\end{equation}
where $F$ is the Faraday constant, $m_{H}$ is the consumed mass of hydrogen and $M_{H}$ is the molar mass of hydrogen. 
The electric current provided by the fuel cell is $I=-\dot Q$. 
The minus sign is due to the fact that the charge $Q$ decreases as the drawn current increases.
The fuel cell output voltage depends on the electrical current and is given by \cite{Larminie}
\begin{equation}\label{operating_voltage}
U_c = E_{oc}-Ir-E_{trans}-E_{act}
\end{equation}
where $E_{oc}$ is the open circuit voltage, $r$ is the internal resistance, $E_{trans}$ corresponds to the mass transport (or concentration) losses, and $E_{act}$ represents activation losses. 
\begin{figure}
	\centering
	\includegraphics[width=0.8\linewidth]{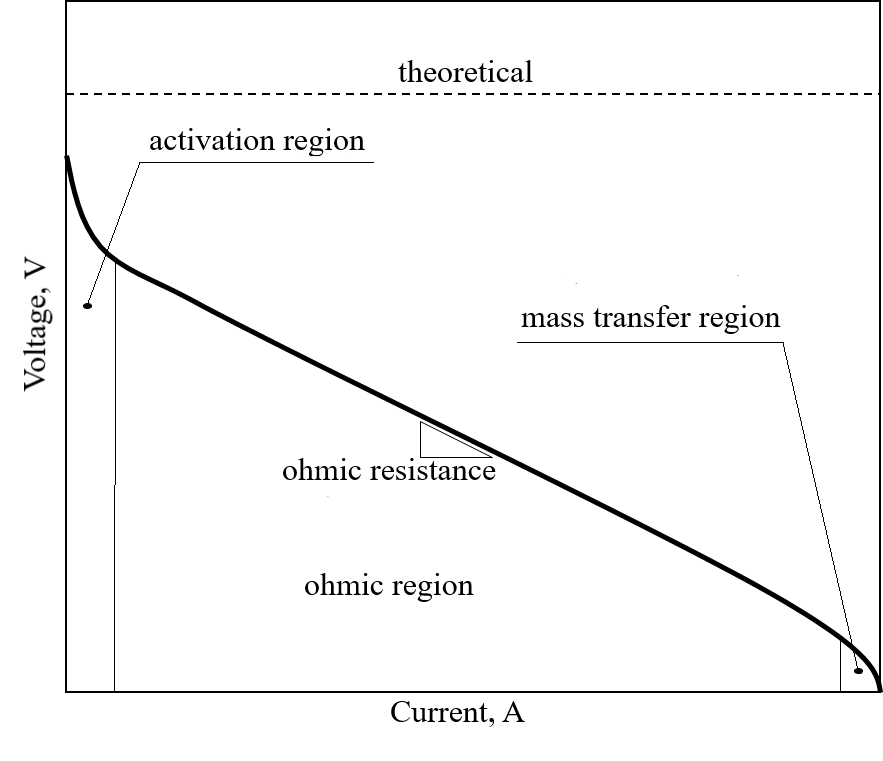}
	\caption{Polarization curve of the fuel cell}
	\label{fig:screenshot006}
\end{figure}

\noindent{\bf Fuel Cell Assumptions:}
\begin{enumerate}
\item The mass transport and activation losses can be neglected, i.e., the fuel cell voltage is an affine function of $I$ given by
\begin{equation}\label{open_circuit}
U_c = E_{oc} - Ir
\end{equation}
\item The internal resistance $r$ is small enough so that
\begin{equation}\label{internalr}
U_c>>Ir
\end{equation}
\end{enumerate}
These assumptions are reasonable for the so-called ohmic region of a typical polarization curve of a fuel cell shown in Figure~2.
As seen in the figure, the ohmic region represents most of the voltage-current characteristic. The activation region is valid for very small currents. The mass transfer region is valid for the highest values of current.
In the ideal situation of no internal resistance (formed by electrolyte resistance to ion flow plus ohmic circuit resistance) one would obtain a horizontal line representing a constant voltage, shown as a dashed line in Figure~2.

For small hydrogen-powered aircraft, the thrust can be generated by the propeller engine. 
The net power required to overcome drag is a function of thrust and velocity and is written as \cite{Anderson}
\begin{equation}\label{power_req}
P_R=Tv=\eta P_e
\end{equation}
where $\eta$ is the system efficiency and $P_e$ is the electric power.
If the voltage of each cell in the stack is $U_c$, then the electric power produced by $n$ cells is
\begin{equation} \label{EPower}
P_e=nU_cI
\end{equation}
Combining (\ref{power_req}) and (\ref{EPower}), and recalling that $T=D$ in steady cruise flight,
\begin{equation}\label{Power_chain}
Dv=\eta n U_c I
\end{equation}
Taking (\ref{open_circuit}) into account, equation (\ref{Power_chain}) becomes
\begin{equation}\label{current_quadratic}
rI^2 - E_{oc}I +\frac{Dv}{\eta n} = 0
\end{equation}
Solving (\ref{current_quadratic}) for $I$ yields
\begin{equation}\label{current_root}
I = \frac{E_{oc} - \sqrt{E_{oc}^2 - 4rDv/(\eta n)}}{2r}
\end{equation}
Note that only the smallest positive real root has been considered because transmitting electricity at a lower curent reduces the fraction of energy lost in the internal resistance through Joule's effect.
Therefore, replacing the drag $D$ from (\ref{Drag}) into (\ref{current_root}) yields an equation for the electric current as an explicit function of $v$ given by
\begin{equation}
I = \frac{E_{oc}}{2r}-\frac{\sqrt{(\eta n E_{oc})^2-2\eta n r\left(C_{D,0}\rho Sv^3+\frac{4KW^2}{\rho Sv}\right)}}{2\eta nr}
\end{equation}
Recalling that $I=-\dot Q$ from (\ref{EPower}), the time rate of change of the charge is
\begin{equation}\label{Qdot}
\dot Q = - \frac{E_{oc}}{2r}+\frac{\sqrt{(\eta n E_{oc})^2-2\eta n r\left(C_{D,0}\rho Sv^3+\frac{4KW^2}{\rho Sv}\right)}}{2\eta nr}
\end{equation} 

For each cell, taking the time derivative of each side of equation~(\ref{single_cell_charge}) and rearranging terms we get the time rate of change of mass of hydrogen fuel
\begin{equation}\label{mass_rate}
\dot{m}_{H} = - \frac{IM_{H}}{2F}
\end{equation}
where $I=-\dot Q$ is the electric current. 
Replacing (\ref{current_root}) in (\ref{mass_rate}) yields (drag is kept as $D$ to make the expression smaller)
\begin{equation}
\dot{m}_{H} = -\frac{M_{H}\left[E_{oc} - \sqrt{E_{oc}^2 - 4Drv/(\eta n)}\right]   }   {4rF}
\end{equation}
The time rate of change of aircraft weight must consider the mass of hydrogen consumed in all $n$ cells in the stack and is then given by 
\begin{equation}\label{weight}
\dot{W}=-n\dot{m}_{H}g = \frac{n M_{H}g\left[E_{oc} - \sqrt{E_{oc}^2 - 4Drv/(\eta n)}\right]   }   {4rF}
\end{equation}
We can now proceed to formulate the Optimal Control Problem (OCP) for steady cruise flight of a hydrogen-powered aircraft.

\section{Problem Formulation}

The OCP is concerned with minimizing the Direct Operating Cost (DOC).
This cost performs a trade-off between energy costs and flight time costs.
The cost of energy can be computed either relative to the consumed hydrogen or relative to the consumed electric charge.
The DOC using the consumed hydrogen is
\begin{equation} \label{DOC}
DOC=\int_0^{t_f}\left(C_t+C_{H}\dot{W} \right) dt
\end{equation}
where $C_H>0$ and $C_t>0$ are the cost of hydrogen fuel per unit weight and time-related costs, respectively.
From equation (\ref{single_cell_charge}), the consumed mass of hydrogen is directly proportional to the charge.
Therefore, the DOC based on the cost of charge can also be written as (\ref{DOC}) but using a different parameter $C_{H}$ to include the constant of proportionality between hydrogen mass and charge.
Since $C_H$ is positive, minimizing $DOC$ is equivalent to minimizing the cost
\begin{equation} \label{Functional}
J=\int_0^{t_f}\left(C_I+\dot W\right) dt
\end{equation}
where $C_I=C_t/C_H$ is the so-called cost index and acts as a trade-off parameter.

The cost index provides an effective and flexible tool to bias the solution between the minimum energy consumption mode and the minimum flight time mode. Speed schedules for optimal operation of the aircraft in flight are determined as a function of $C_I$, so that for any value of $C_I$ chosen by the pilot one can obtain the most economical speed.
A larger $C_I$ corresponds to a higher true airspeed but results in the accelerated consumption of hydrogen (or electric charge). 
Combining the cost functional (\ref{Functional}) with the flight dynamics (\ref{Drag})--(\ref{Dynamics_reduced}) and (\ref{weight}) leads to the following OCP:
\begin{equation} \label{OCP}
\begin{split}
J^*=&\inf_{v, t_f}\int_0^{t_f}(\dot{W}+C_I) dt\\
&s.t.\\
&\dot{x}=v\\
&\dot{W} = \frac{n M_{H}g\left[E_{oc} - \sqrt{E_{oc}^2 - 4Drv/(\eta n)}\right]   }   {4rF} \\
&D=\frac{1}{2}C_{D,0}\rho Sv^2+\frac{2KW^2}{\rho Sv^2}\\
&x(0)=0, \ x(t_f)=x_d,\ W(0)=W_0,\ v>0\\
\end{split}
\end{equation}
where the final time $t_f$ is an optimization variable.
This problem formulation fits the general setting of a Lagrange OCP given by \cite{Liberzon}
\begin{eqnarray}\label{OCPLagrange}
&\min\limits_{u(t),t_f} & \int_{0}^{t_f} L(X(t),u(t)) dt\nonumber\\
&\mbox{s.t.}~&\dot X(t) =f(X(t),u(t))\nonumber\\
&&~~~~0=\psi(X(0),X(t_f))
\end{eqnarray}
with the additional input constraint $u(t)=v(t)>0$. 
Since the constraint corresponds to a strict inequality the solution of the problem can be carried without adding the constraint and using it only at the end to validate the candidate solutions.
In (\ref{OCP}) we have chosen the state $X(t)=[x(t)~~W(t)]^T$ and
\begin{eqnarray}\label{substitutions}
L(X(t),u(t))  &=& \dot{W}+C_I\nonumber\\
f(X(t),u(t)) &=& \left [v\quad\frac{n M_{H}g\left[E_{oc} - \sqrt{E_{oc}^2 - 4Drv/(\eta n)}\right]} {4rF}\right]^T\nonumber\\
\psi(X(0),X(t_f)) &=& \left[x(0)\quad x(t_f)-x_d\quad W(0)-W_0\right]^T
\end{eqnarray}

\section{Problem Solution}

\begin{theorem}\label{HydroFMS}
	Under the assumptions (\ref{open_circuit})--(\ref{internalr}), the optimal solution $v^*$ of the OCP stated in (\ref{OCP}) for steady cruise flight at a constant altitude is a positive root (if it exists) of the nonlinear equation 
	\begin{eqnarray}\label{solution}
	(1+J_W)\frac{M_Hg}{2F}\left[\frac{n\left(E_{oc} - \sqrt{E_{oc}^2 - 4Drv/(\eta n)}\right)}{2r} -\frac{D_vv^2+Dv}{\eta \sqrt{E_{oc}^2- 4Drv/(\eta n)}}\right]+C_I =0&\nonumber\\
	\end{eqnarray}
	where the drag $D$ is given by (\ref{Drag}).
\end{theorem}
\vspace{10pt}

{\bf Proof:}
For a problem in the form (\ref{OCPLagrange}) the Hamiltonian is defined to be
\begin{equation}\label{theHamiltonian}
H(X(t),u(t),\lambda(t))=L(X(t),u(t))+\lambda^T(t)f(X(t),u(t)).
\end{equation}
where the costate for problem (\ref{OCP}) is $\lambda(t)=[J_x(t)~~J_W(t)]^T$.
In the rest of the proof we may ommit the direct dependence of $X, \lambda,$ and $u$ on time for simplicity.
From (\ref{theHamiltonian}) and (\ref{substitutions}), the Hamiltonian can be written as
\begin{equation} \label{hamiltonian}
\begin{split}
&H (x,W,v,J_x,J_W)= (1+J_W)\frac{n M_{H}g\left[E_{oc} - \sqrt{E_{oc}^2 - 4Drv/(\eta n)}\right]   }   {4rF} \\
&+J_xv + C_I
\end{split}
\end{equation}
Since the Hamiltonian is explicitly independent of time and the terminal time is free, then $\dot H=0$ and from the transversality conditions \cite{Pontryagin,Liberzon}
\begin{equation} \label{PMP}
H(x(t),W(t),v^*(t),J_x(t),J_W(t))=H(x(t_f),W(t_f),v^*(t_f),J_x(t_f),J_W(t_f))=0
\end{equation}
Additionaly, the costate (Hamilton's) equations must be satisfied  \cite{Pontryagin,Liberzon} and can be written for a continuously differentiable Hamiltonian $H$ as
\begin{equation} \label{hamilton_eqs}
\begin{split}
&\frac {dJ_x}{dt}=-\frac {\partial {H}} {\partial {x}}=0\\
&\frac {dJ_W}{dt}=-\frac {\partial {H}} {\partial {W}}
\end{split}
\end{equation}
Combining equations (\ref{hamiltonian}) and (\ref{PMP}) yields
\begin{equation} \label{H=0}
\begin{split}
&(1+J_W)\frac{n M_Hg\left[E_{oc} - \sqrt{E_{oc}^2 - 4Drv/(\eta n)}\right]   }   {4rF} \\
&+J_xv + C_I = 0
\end{split}
\end{equation}
From (\ref{hamilton_eqs}), one can oberve that $J_x$ must be a constant. 
According to Pontryagin's minimum principle \cite{Pontryagin, Liberzon}, the necessary condition for a minimum is
\begin{equation} \label{NC}
\frac{\partial{H}}{\partial{v}} = (1+J_W)\frac{M_{H}g (D_vv+D)    }   {2\eta F\sqrt{E_{oc}^2 - 4Drv/(\eta n)}}+J_x = 0
\end{equation}
where the notation $D_v$ refers to the partial derivative of the drag $D$ relative to the velocity $v$.
Note that the denominator in (\ref{NC}) is equal to $2\eta F(U_c-rI)$ from equations (\ref{current_root}) and (\ref{open_circuit}). Moreover, the denominator is strictly positive from assumption (\ref{internalr}).
Therefore, the partial derivative (\ref{NC}) is well defined.
Solving (\ref{NC}) for $J_x$ leads to
\begin{equation} \label{Jx}
J_x= - (1+J_W)\frac{M_Hg (D_vv+D)    }   {2\eta F\sqrt{E_{oc}^2 - 4Drv/(\eta n)}}
\end{equation}
Replacing (\ref{Jx}) in (\ref{H=0}) and factoring common terms yields (\ref{solution}).
\done
\vspace{10pt}

\textbf{Remark 1}: The solutions of (\ref{solution}) can be obtained from the positive real roots of an eighth order polynomial in $v$ after multiplying (\ref{solution}) by the square root term, some rearranging, and squaring both sides of the resulting equation.

\textbf{Remark 2} An important difference between hydrogen-powered aircraft and all-electric battery-powered aircraft is that the costate $J_W$ does not appear in the solution for all-electric aircraft because the weight of battery-powered aircraft is constant. This is not the case for hydrogen-powered aircraft because hydrogen is being consumed. Note that the costate $J_W$ is the sensitivity of the optimal cost-to-go with respect to aircraft weight. In reference\cite{Sorensen}  the costate $J_W$ is assumed to be zero because the cost-to-go does not change much with the aircraft weight. In fact, for aircraft powered by fuel cells the total aircraft weight is considerably greater than that of the hydrogen fuel stored onboard. This assumption is valid for the prototype HY4 aircraft used in the simulations, which has a maximum takeoff weight around $1500~kg$ and for which the hydrogen weight component is $9~kg$ \cite{HY4}.

\begin{figure}[t]
	\centering
	\includegraphics[width=0.9\linewidth]{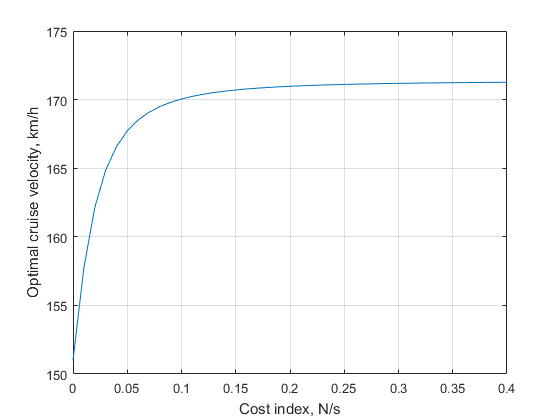}
	\caption{Velocity as a function of cost index.}
	\label{fig:optvelo}
\end{figure}
\vspace{10pt}

\section{Simulation Results}
The aircraft model of the H4Y hydrogen-powered aircraft will be used for simulations in Matlab/Simulink. 
The HY4’s drive train consists of a hydrogen storage unit, a low-temperature hydrogen fuel cell and a high performance battery. The fuel cell converts the energy of the hydrogen directly into electrical energy. The only waste product in the process is clean water. The electric motor uses the power generated to propel the aircraft. If the hydrogen required for the fuel cell is generated via electrolysis using power from a renewable energy source, the HY4 can fly without generating any emissions at all.
The technical data is taken from \cite{HY4} and is presented in Table I. 
The initial and final conditions are shown in Table \ref{ICFC}.
Bearing in mind the comment on remark 1, a suboptimal solution to the OCP can be derived by assuming $J_W=0$ and by solving the corresponding equation that appears in theorem \ref{HydroFMS} after transforming it in a polynomial root finding problem.
The main advantage of the suboptimal solution is that it allows one to compute the true airspeed using Newton's method, whereas the optimal solution requires more sophisticated techniques such as the shooting method. 

\begin{figure}[t]
	\centering
	\includegraphics[width=0.9\linewidth]{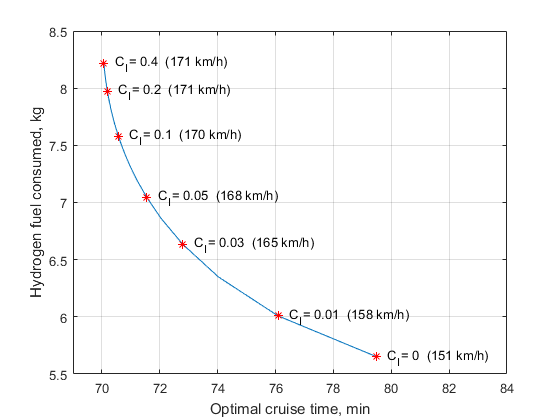}
	\caption{Pareto-optimal trade-off curve.}
	\label{fig:pareto}
\end{figure}
\begin{table}[h]
	\begin{center}
		\caption{HY4 Data}\label{Technical Data}
		\begin{tabular}{ll}
			\hline
			& TECHNICAL DATA\\
			\hline
			Wingspread & 21.36 $m$ \\
			Length & 7.4 $m$ \\
			Maximum takeoff weight & 1500 $kg$ \\
			Maximum speed & 200 $km/h$ \\
			Cruise speed & 145 $km/h$ \\
			Range & 800--1500 $km$ \\
			Motor type & Permanent magnet synchronous \\
			Motor power & 80 $kW$ \\
			Total system efficiency & 0.44 \cite{Hepperle} \\
			Parasitic drag coefficient & 0.025 (estimate)\\
			Induced drag constant & 0.039 (estimate)\\
			\hline	
			& HYDROGEN STORAGE SYSTEM \\
			\hline
			Pressure & 437 $bar$ \\
			Hydrogen amount & 9 $kg$ \\
			Weight & 170 $kg$ \\
			\hline
			& FUEL CELLS DATA \\	
			\hline
			Type & Low temperature PEM \\
			Nominal power & 45 $kW$ \\
			Quantity cells & 440 \\
			Quantity modules & 4 \\
			Cathode gas & Air \\
			Total weight & 100 $kg$ \\
			Ohmic resistance & 0.005 $Ohm$ (estimate) \\
			Open circuit voltage & 1.1 $V$ (estimate) \\
			\hline
			& BATTERY DATA \\	
			\hline
			Capacity & 21 $kWh$ \\
			Quantity cells & 80 (4$\times$ 20) \\
			Type & Li-Po \\
			Nominal power & 45 $kW$ \\
			Maximum current & 600 $A$ \\
			Modules & 4 \\
			Weight & 130 $kg$\\
			\hline		
		\end{tabular}
	\end{center}
\end{table}

\begin{table}[h]
	\begin{center}
		\caption{Simulation initial and final conditions}\label{ICFC}
		\begin{tabular}{cc}
			\hline
			Cruising altitude & 1 $km$ \\
			Final position &  200 $km$ \\
			Hydrogen fuel weight &  9 $kg$ \\
			Weight &  1500 $kg$ \\  
			\hline
		\end{tabular}
	\end{center}
\end{table}

The graph shown in Figure~\ref{fig:optvelo} gives the optimal cruise velocity as a function of cost index. It is clear from the plot that, as expected, larger values of $C_I$ correspond to a higher cruising velocity.
The Pareto-optimal trade-off curve between the final optimal time and the amount of hydrogen fuel consumed is shown in Figure~\ref{fig:pareto}. From this figure we can see that increasing the speed from $151~km/h$ to $171~km/h$ only increases the consumption of hydrogen fuel by less than $5kg$ while saving roughly $20$ minutes of flight. At a velocity of $171~km/h$ the derivative of the trade-off curve becomes much steeper and therefore the amount of time of flight saved by increasing speed beyond $171~km/h$ will be significantly reduced.

\section{Conclusions}
This paper has formulated and solved a direct operating cost minimization problem for hydrogen-powered aircraft in steady cruise at a constant altitude. The solution for the optimal speed corresponds to the optimal trade-off between energy consumption and flight time.  
A case study using the HY4 airplane model parameters has illustrated the methodology.
The results show that an airplane powered by hydrogen fuel cells is a potential solution for the emerging business of urban mobility.
One huge potential of this solution is that the only exhaust products generated by hydrogen-powered aircraft are heat and water vapor. 
Additionally, a fuel cell aircraft produces zero noise, which is especially important during takeoff and landing in urban centers.
The proposed solution provides a Pareto-optimal trade-off curve that helps in the selection of the airspeed that yields a good compromise between energy consumption and time of flight.

\end{document}